\documentclass[a4paper,11pt,twoside]{amsart}
\usepackage[latin1]{inputenc}
\usepackage[T1]{fontenc}
\usepackage{amsmath}
\usepackage{amsthm}
\usepackage{amssymb}
\usepackage[all]{xy}

\newtheorem{thm}{Theorem}[section]

\newtheorem{prop}[thm]{Proposition}

\newtheorem{lem}[thm]{Lemma}
\newtheorem{cor}[thm]{Corollary}

\numberwithin{equation}{section}

\theoremstyle{definition}

\newtheorem{remark}[thm]{Remark}

\newcommand{\qqed}{\hspace*{\fill}$\Box$}

\newcommand{\Db}{{\rm D}^{\rm b}}

\newcommand{\Aut}{{\rm Aut}}

\newcommand{\NS}{{\rm NS}}
\newcommand{\Pic}{{\rm Pic}}

\newcommand{\coh}{{\rm Coh}}

\newcommand{\End}{{\rm End}}
\newcommand{\Hom}{{\rm Hom}}
\newcommand{\Spec}{{\rm Spec}}

\newcommand{\Ext}{{\rm Ext}}
\newcommand{\CH}{{\rm CH}}

\newcommand{\cal}{\mathcal}
\newcommand{\ka}{{\cal A}}

\newcommand{\kc}{{\cal C}}
\newcommand{\kd}{{\cal D}}

\newcommand{\ko}{{\cal O}}

\newcommand{\ks}{{\cal S}}

\newcommand{\ZZ}{\mathbb{Z}}
\newcommand{\QQ}{\mathbb{Q}}

\newcommand{\CC}{\mathbb{C}}

\renewcommand{\to}{\xymatrix@1@=15pt{\ar[r]&}}
\renewcommand{\rightarrow}{\xymatrix@1@=15pt{\ar[r]&}}
\renewcommand{\mapsto}{\xymatrix@1@=15pt{\ar@{|->}[r]&}}
\renewcommand{\twoheadrightarrow}{\xymatrix@1@=15pt{\ar@{->>}[r]&}}
\renewcommand{\hookrightarrow}{\xymatrix@1@=15pt{\ar@{^(->}[r]&}}
\newcommand{\congpf}{\xymatrix@1@=15pt{\ar[r]^-\sim&}}
\renewcommand{\cong}{\simeq}

\begin{document}
\title{A note on the Bloch--Beilinson conjecture for K3 surfaces and spherical objects}
\author{D. Huybrechts}
\address{Mathematisches Institut,
Universit{\"a}t Bonn, Endenicher Allee 60, 53115 Bonn, Germany}
\email{huybrech@math.uni-bonn.de} 
\thanks{This work was supported by the SFB/TR 45 `Periods,
Moduli Spaces and Arithmetic of Algebraic Varieties' of the DFG
(German Research Foundation). The hospitality and financial support
of the Mathematical Institute Oxford is gratefully acknowledged.
}
\maketitle

For a projective K3 surface $X$ over an algebraically closed field $k$ let $\Db(X)$ denote the bounded
derived category of coherent sheaves. Spherical objects, e.g.\ line bundles and rigid stable bundles,
play a central role in the study of $\Db(X)$ and, as it turns out, also of $\CH^*(X)$. Our aim here
is to make this more precise by studying the $k$-linear triangulated subcategory $\ks^*\subset\Db(X)$
generated by all spherical objects (see below for definitions).
The main result concerns K3 surfaces over number fields and is a consequence of Thomason's classification of dense subcategories. It can be stated as follows.

\begin{thm}\label{thm:MAIN}
Let $X$ be a K3 surface over $\bar\QQ$ with Picard number $\rho(X)\geq2$. Then the following conditions are
equivalent:\\
i) $\CH^2(X)\cong\ZZ$.\\
ii) $\Db(X)=\ks^*$, i.e.\ $\Db(X)$ is generated by spherical objects.\\
iii) The triangulated category $\ks^*$ admits a bounded t-structure. 
\end{thm}

Note that i) is predicted by the Bloch--Beilinson conjecture applied to K3 surfaces. The implication
ii) $\Rightarrow$ iii) is obvious and ii) $\Rightarrow$ i) follows from \cite{HuyJEMS}. This note is concerned with the
curious observation that i) and ii) are in fact equivalent (and that both are implied by iii)).

In \cite{HuyJEMS,HuyMSRI} and elsewhere we have alluded to ii) as a `categorial
logical possibility', in analogy to Bogomolov's `logical possibility' that every
$\bar\QQ$-rational point on a K3 surface might be contained in a rational curve (see \cite{Bog}).
Whether either of the two is a way towards a proof of
the Bloch--Beilinson conjecture for K3 surfaces seemed
doubtful, but due to Theorem \ref{thm:MAIN} one now knows at least that our categorial one is strictly
equivalent to it.

Theorem \ref{thm:MAIN} is valid for any algebraically closed field $k$ of characteristic zero, but for
$k=\CC$, and in fact whenever ${\rm trdeg}( k)>0$,
a result of Mumford shows that none of the conditions, i)-iii) holds.

For a discussion
of spherical objects, Chow groups and derived categories we refer to \cite{HuyMSRI}. The notion of spherical objects and the precise definition of $\ks^*$ can be found in Section \ref{sec:SECT1}. The equivalence of ii) and iii) is proved in Section \ref{sec:SECT2}. The case $k=\bar\QQ$ is discussed in Section \ref{sect:Qbar}.
  
\section{The spherical category $\ks^*$}\label{sec:SECT1}

Let $X$ be a K3 surface over a field $k$. Then the category $\coh(X)$ will be considered
as a $k$-linear abelian category and its bounded derived category
$\Db(X):=\Db(\coh(X))$ as a $k$-linear triangulated category.

Recall that an object $A\in\Db(X)$ is called \emph{spherical} if there
exists a $k$-linear isomorphism $\Ext^*(A,A)\cong H^*(S^2,k)$. In other words
$\Ext^i(A,A)$ is one-dimensional for $i=0,2$ and zero otherwise.
We denote the collection of all spherical objects by $$\ks\subset\Db(X).$$

We can view $\ks$ simply as a set or as a full subcategory of $\Db(X)$. Note that
$\ks$ is invariant under shift, but it is not a triangulated subcategory of $\Db(X)$
and, in fact, not even an additive one. Indeed, the direct sum $A_1\oplus A_2$ of two
spherical objects $A_1,A_2\in\ks$ is clearly not spherical. Instead we consider the category
$\ks^*$ generated by $\ks$.

To avoid confusion we shall spell out the definition of $\ks^*$. It is the full subcategory of $\Db(X)$
with objects $\bigcup\ks^{n}$, where $\ks^{n}$ is defined recursively by $\ks^{1}:=\ks$ and
$\ks^{n+1}:=\ks^n\ast\ks$. Here for two subcategories $\ka_1,\ka_2\subset\Db(X)$ one
lets $\ka_1\ast\ka_2$ be the full subcategory of all objects $A\in\Db(X)$ such that there exist
objects $A_i\in\ka_i$, $i=1,2$, and an exact triangle $A_1\to A\to A_2$. For completeness sake we 
give a proof of the rather obvious

\begin{lem}
$\ks^*\subset\Db(X)$ is a full triangulated subcategory.
\end{lem}

\begin{proof}
It is enough to show that for objects $A\in\ks^m,B\in\ks^n$ and an exact triangle
$A\to C\to B$ one has $C\in\ks^{m+n}$. This can be proved by induction on $n$. If $n=1$, i.e.\
$B\in\ks$, then $C\in\ks^{m+1}$ by definition. Otherwise, there exists an exact triangle
$B_1\to B\to B_2$ with $B_1\in\ks^{m-1}$ and $B_2\in\ks$. Then consider the
diagram
$$\xymatrix{A\ar@{=}[d]\ar[r]&C_1\ar[d]\ar[r]&B_1\ar[d]\\
A\ar[r]&C\ar[r]\ar[d]&B\ar[d]\\
&B_2\ar@{=}[r]&B_2}$$
with exact triangles $A\to C_1\to B_1$ and $C_1\to C\to B_2$. From the first one and the induction hypothesis
one deduces $C_1\in\ks^{m+n-1}$. The second one and the definition of $\ks^{m+n}$ then yield the assertion.
\end{proof}

\begin{remark}
Note that $\ks$ and $\ks^*$ are \emph{strictly full}, i.e.\ any object in $\Db(X)$ isomorphic to an object in $\ks$ or $\ks^*$  is contained in $\ks$ resp.\ $\ks^*$. 

The category $\ks^*$ can equivalently be described as the smallest strictly full triangulated subcategory of $\Db(X)$ that contains $\ks$. But note that we do not require $\ks^*$ to be closed under taking direct summands, i.e.\ if $E_1,E_2\in\Db(X)$ with $E_1\oplus E_2\in\ks^*$, then $E_1,E_2$ need not necessarily be objects in $\ks^*$.

There are various notions of generators of a triangulated category. E.g.\ one also says
that a collection of objects $\{E_i\}$ generates $\Db(X)$ if $\Ext^n(E_i,E)=0$ for all $i,n$ implies
$E\cong 0$. E.g.\ for an ample line bundle $\ko(1)$ on $X$ the set $\{\ko(i)\}$ generates
$\Db(X)$ in this sense (cf.\ Lemma \ref{lem:dense1}) and hence also $\ks^*$. But the notion used here is 
different. 
\end{remark}

It is easy to see that $\ks^*$ is \emph{dense} in  $\Db(X)$, i.e.\ that for every object $E\in\Db(X)$ there exists an object $F\in\Db(X)$ with $E\oplus F\in\ks^*$. This is equivalent to 
$\Db(X)$ being the smallest strictly full thick triangulated subcategory of $\Db(X)$ that contains $\ks$. The argument to prove that $\ks^*$ is dense can be split in two steps and is well-known (see e.g.\ the more general  \cite[Lem.\ 3.12]{Thom}). Again for completeness sake we include the details.

\begin{lem}\label{lem:dense1}
Let $X$ be a smooth projective variety  of dimension $d$ with an ample line bundle $\ko(1)$.
Then for any $E\in\coh(X)$ there exists $F\in\coh(X)$ and a complex $$K^\bullet=[\ko(-n_d)^{\oplus m_d}\to\ldots\to
\ko(-n_0)^{\oplus m_0}]$$ such that $$E\oplus F[d]\cong K^\bullet.$$
\end{lem}

\begin{proof}
Choose a resolution of $E$ of the form $$
\ldots\to\ko(-n_1)^{\oplus m_1}\to\ko(-n_0)^{\oplus m_0}\to E\to0$$
and let $F$ be the kernel of $\ko(-n_d)^{\oplus m_d}\to\ko(-n_{d-1})^{\oplus m_{d-1}}$. The extension class
of the exact triangle $F[d]\to K^\bullet\to E$ is an element in $\Hom(E,F[d+1])=\Ext^{d+1}(E,F)=0$. Thus,
$K^\bullet\cong E\oplus F[d]$.
\end{proof}

\begin{lem}\label{lem:dense2} 
Let $\kd:=\Db(\kc)$ be the bounded derived category of an abelian category $\kc$. A strictly full triangulated
subcategory $\kd_0\subset\kd$ is dense if and only if for all $E\in\kc$ there exists an object $F\in\kd$ with
$E\oplus F\in\kd_0$.
\end{lem}

\begin{proof} The `only if' is clear. Let us prove the `if'.
For any $E\in \kd$ we have to find a complex $F\in\kd$ such that $E\oplus F\in\kd_0$ which we will do by induction
on the length $\ell(E)$ of the complex $E$. Choose an exact triangle $E_1\to E\to E_2$ with $E_2\in\kc[n]$ for some $n$ and such that $\ell(E_1)<\ell(E)$. By induction hypothesis there exist $F_1,F_2\in\kd$
such that $E_i\oplus F_i\in\kd_0$, $i=1,2$. Taking the direct sum of the exact triangles
$E_1\to E\to E_2$, $F_1=F_1\to 0$, and $0\to F_2=F_2$ yields an exact triangle of the form
$E_1\oplus F_1\to E\oplus F_1\oplus F_2\to E_2\oplus F_2$. Since $E_i\oplus F_i\in\kd_0$ and $\kd_0$ is
closed under extensions, this proves $E\oplus F_1\oplus F_2\in\kd_0$. 
\end{proof}

\begin{cor}
$\ks^*\subset\Db(X)$ is a dense strictly full triangulated subcategory.
\end{cor}

\begin{proof}
This follows directly from the above using the fact that all line bundles, and in particular the line bundles $\ko(-n_i)$,
are spherical. Thus, the complex $K^\bullet$ occurring in Lemma \ref{lem:dense1} is an object in $\ks^*$. 
\end{proof}

Next, let us recall Thomason's classification of dense subcategories. Here and in the sequel
$K(\kd)$ of a triangulated or an abelian category $\kd$ will be the Grothendieck group of $\kd$. If $\kd=\Db(\kc)$
with $\kc$ abelian, then $K(\kd)\cong K(\kc)$.

\begin{thm}{\bf (Thomason, \cite{Thom})} Let $\kd$ be an essentially small triangulated category. Then there exists a natural
bijection 
$$\{\text{strictly full dense triangulated subcategories}\}\leftrightarrow\{\text{subgroups } H\subset K(\kd)\}.$$ 
\end{thm}

The bijection is given by mapping a triangulated subcategory $\kd_0\subset\kd$ to the image of
the natural map $K(\kd_0)\to K(\kd)$. The inverse map sends a subgroup $H\subset K(\kd)$ to
the strictly full dense triangulated subcategory $\kd_H\subset \kd$ of all objects $E\in\kd$ with
$[E]\in H$.

Moreover, for a dense triangulated subcategory $\kd_0\subset\kd$ the map $K(\kd_0)\to K(\kd)$ is in fact
injective, which was proved by Thomason as a consequence of the above theorem (cf.\ \cite[Cor.\ 2.3]{Thom}).
Applied to our situation one finds that the natural $K(\ks^*)\to K(X)\cong K(\Db(X))$ identifies
$K(\ks^*)$ with a subgroup of $K(X)$. Furthermore, the theorem yields the following characterization
of $\ks^*$.

\begin{cor}\label{cor:sphcatThom} The subcategory
$\ks^*\subset\Db(X)$ is the strictly full dense subcategory of all objects $E\in\Db(X)$ with $[E]\in K(\ks^*)\subset K(X)$.\qqed
\end{cor}

In particular, for arbitrary $E\in\Db(X)$ the object $E\oplus E[1]$ is contained in $\ks^*$ which seems difficult
to prove directly.

Let us recall the following result from \cite{HuyJEMS}.

\begin{thm}\label{thm:JEMS}
If $k$ is algebraically closed of characteristic zero and $\rho(X)\geq2$, then $K(\ks^*)=R(X)$.
\end{thm}

Here, $R(X)$ is the Beauville--Voisin ring generated by the Chern characters ${\rm ch}(L)$ of all line bundles
$L\in\Pic(X)$ and the distinguished class $c_X\in\CH^2(X)$ satisfying $24c_X=c_2(X)$.\footnote{Note that here and in the sequel we implicitly identify
the Chow ring $\CH^*(X)$ with the Grothendieck group $K(X)$, i.e.\ $\CH^*(X)=K(X)$,  tensoring with $\QQ$ 
is not needed for K3 surfaces.}
Alternatively, $R(X)$ is spanned
by the Mukai vectors $v(L)\in\CH^*(X)$ of line bundles and the theorem says that the Mukai vector $v(E)$ of
all spherical objects $E$ are also contained in $R(X)$. As was shown in \cite{BV}, the cycle map induces an
isomorphism $R(X)\cong \ZZ\oplus\NS(X)\oplus\ZZ$. 

Thus under the assumption of the theorem Corollary \ref{cor:sphcatThom} says that $$\ks^*=\kd_{R(X)}.$$ Of course the inclusion $\ks^*\subset \kd_{R(X)}$ follows immediately from Theorem \ref{thm:JEMS}, but for the equality Thomason's result is used.

\section{$\Db(X)$ versus $\ks^*$}\label{sec:SECT2}

Here we shall discuss some of the similarities between
$\ks^*$ and $\Db(X)$ and some of their differences. As we will explain in 
Section \ref{sect:Qbar}, one expects $\ks^*=\Db(X)$  for $k=\bar\QQ$ 
but for $k=\CC$ the two categories are very different. But, nevertheless, even for $k=\CC$ 
dealing with $\Db(X)$ or $\ks^*$ seems equivalent for many purposes.


The main difference between $\Db(X)$ and $\ks^*$ is that $\Db(X)$, as the bounded derived category
of an abelian category, is Karoubian. This not obvious and in general not true for $\ks^*$.
Recall that a triangulated category $\kd$ is \emph{Karoubian} if it is idempotent-split, i.e.\ every morphism $f:E\to E$ with $f^2=f$ comes from a direct sum decomposition $E={\ker}(f)\oplus{\rm Im}(f)$.  Thus, since
$\ks^*\subset\Db(X)$ is dense, the triangulated category $\ks^*$ is Karoubian if and only if $\ks^*=\Db(X)$.
In general, $\Db(X)$ can be seen as the Karoubian (or idempotent) closure of $\ks^*$, i.e.\
$\Db(X)$ is equivalent to the category of all pairs $(E,f)$ with $E\in\ks^*$ and an idempotent
$f\in\End(E)$ (see e.g.\ \cite{BalmSchl}).

\subsection{} Let us first study autoequivalences of $\Db(X)$ and of $\ks^*$.
Consider  two projective K3 surfaces  $X$ and $Y$ 
over an algebraically closed field $k$. Let $\ks_X^*\subset\Db(X)$ and $\ks^*_Y\subset\Db(Y)$ be the corresponding
subcategories generated by the collections $\ks_X$ resp.\ $\ks_Y$ of
spherical objects. Let now $\Phi:\Db(X)\congpf\Db(Y)$ be a linear
exact equivalence. Since the notion of a spherical object is purely categorial, $\Phi$ induces a linear equivalence
$\Phi_\ks:\ks_X\congpf\ks_Y$ and an exact linear equivalence $\Phi_{\ks^*}:\ks_X^*\congpf\ks_Y^*$.
In particular, this yields a natural homomorphism $$\Aut(\Db(X))\to\Aut(\ks^*_X),$$
where $\Aut(\kd)$ for a linear triangulated category $\kd$ denotes the group
of isomorphism classes of all exact linear autoequivalences.

\begin{prop}
Any exact linear equivalence $\Psi:\ks_X^*\congpf\ks_Y^*$ is isomorphic to $\Phi_{\ks^*}$ for some
uniquely determined  $\Phi:\Db(X)\congpf\Db(Y)$ (up to isomorphism).
\end{prop}

\begin{proof}
Any (exact) equivalence between two categories extends naturally to an (exact)
equi\-valence of their Karoubian closures (see \cite{BalmSchl}).
Since the Karoubian closure of $\ks_X^*$ is $\Db(X)$ and similarly for $Y$, any given
exact equivalence $\Psi:\ks_X^*\congpf\ks_Y^*$ extends naturally to an exact equi\-valence
$\Db(X)\congpf\Db(Y)$. By the same argument or by a result of Bondal and Orlov (see e.g.\ \cite[Ch.\ 4.3]{HuyFM}) 
an equivalence $\Phi:\Db(X)\congpf\Db(Y)$ is uniquely determined by its restriction to $\ks_X^*$.
\end{proof}

This yields in particular the following
\begin{cor} For K3 surfaces $X$ and $Y$ over an algebraically closed field $k$ one has:\\
i) Derived equivalence is determined by the category of spherical objects. More precisely,
 $$\Db(X)\cong\Db(Y)\text{ if and only if }\ks^*_X\cong\ks_Y^*.$$
ii) Restricting an autoequivalence of $\Db(X)$ to $\ks_X^*$ yields an isomorphism
$$\Aut(\Db(X))\congpf\Aut(\ks_X^*).$$ 

\end{cor}

\begin{remark}
Can one replace $\ks^*_X$ by $\ks_X$? 

i) If $\ks_X$ is viewed as a linear category then the same result
of Bondal and Orlov shows that $\Aut(\Db(X))\to\Aut(\ks_X)$ is injective.

ii) If $\ks_X$ is  viewed as a simple set, then one can  show that elements in the kernel of $\Aut(\Db(X))\to\Aut(\ks_X)$ are 
of the form $f^*$ for some automorphism $f$ of $X$ (see \cite{HuyStab}). In fact,
it is expected that the kernel coincides with the finite group of all automorphisms of $X$ that
act trivially on $\NS(X)$.

iii) The surjectivity of $\Aut(\Db(X))\to\Aut(\ks_X)$ for $\ks_X$ as a set (or as a quandle)
seems doubtful, as the $k$-linear structure should enter somehow. In this context it would be interesting to understand the action of ${\rm Gal}(k/k_0)$ on $\ks_X$, where $k_0\subset k$ is the field of definition
of $X$. 

For $\ks_X$ as a linear category the surjectivity seems more likely and would depend on whether
$\ks_X$ determines the structure of $\ks_X^*$, both viewed as abstract categories (and not as subcategories of
$\Db(X)$).
\end{remark}

\subsection{} Let us now turn to stability conditions. Here the category $\ks^*$ seems
most unsuitable, as it can be non-Karoubian (e.g.\ for $k=\CC$).
For this first recall the following basic fact, see \cite{LeChen}.

\begin{prop}\label{prop:LeChen}
If a triangulated category $\kd$ admits a bounded t-structure, then it is Karoubian.
\end{prop}

The next result is in contrast to the case $\Db(X)$ for a K3 surface $X$ over $\bar\QQ$
on which stability conditions can be constructed copying the arguments in \cite{BrK3}.

\begin{cor}\label{cor:no-t}
If the inclusion $\ks^*\subset \Db(X)$ is not an equality, then $\ks^*$ does not admit a bounded
t-structure and, hence, no stability condition.
\end{cor}

\begin{proof}
Recall that a stability condition, e.g.\ in the sense of \cite{Br}, is given by a bounded t-structure together with
a stability function with the HN-property on its heart. But if $\ks^*\ne\Db(X)$, then $\ks^*$ is not Karoubian and
by Proposition \ref{prop:LeChen} does not admit a bounded t-structure.
\end{proof}

 Clearly, the corollary proves that  in Theorem \ref{thm:MAIN} iii) implies ii). The converse of it is obvious.

\begin{remark}
Triangulated K3 categories have been studied in other situations and most successfully in the case
of local K3 categories (see e.g.\ \cite{BKl,Ishi,Th}).
More precisely consider the minimal resolution $\pi:X\to\Spec(k[x,y]^G)$ of
an ADE singularity ($G\subset {\rm Sl}(2)$ a finite group). The usual derived category $\Db(X)$,
which can also be seen as the bounded derived category of the abelian category $\coh^G(\Spec(k[x,y]))$
of $G$-equivariant coherent sheaves on the plane, contains two natural triangulated subcategories
$\kd\subset\hat\kd\subset\Db(X)$. The objects in $\hat\kd$ are complexes with cohomology supported
on the exceptional divisor $E\subset X$ and objects $E\in\kd$ satisfy in addition $R\pi_*E=0$.

The category $\hat\kd$ contains spherical objects $A_0,\ldots,A_n$ that correspond to the vertices
of the extended Dynkin graph of $G$ or, equivalently, to the irreducible representations of $G$.
If one lets $A_0$ be the one corresponding to the trivial representation, then $A_1,\ldots,A_n\in\kd$.

It is known that the Grothendieck group  of $\hat \kd$ (resp.\ $\kd$) is freely generated by the classes
 $[A_i]$. In particular, in both cases the collection of spherical objects $\hat\ks\subset\hat\kd$ (resp.\ $\ks\subset\kd$)
 and its associated category $\hat\ks^*$ (resp.\ $\ks^*$) spans the full Grothendieck group. Thus, by Corollary
 \ref{cor:no-t} $\hat\ks^*=\hat\kd$ (resp.\ $\ks^*=\kd$), i.e.\ both categories are generated
 by spherical objects (not allowing taking direct summands). Presumably, this can be also proved more directly,
 but it is instructive to see how everything falls into place in this geometrically easier case.
 
 Similar arguments would apply to another local K3 category, the derived category of complexes on the total
 space of $\ko(-2)$ which are concentrated in the zero section.
\end{remark}

In \cite{HuyStab} we show that stability conditions on $\Db(X)$ are determined by their behaviour
with respect to spherical objects, so morally by their `restriction' to $\ks^*$. In order to make this precise
we introduce a modified notion of stability conditions that is applicable to $\ks^*$.  
\section{Special base fields}\label{sect:Qbar}

Due  to a result of Mumford \cite{Mum} one knows that for $k=\CC$ the Chow group $\CH^2(X)$ is infinite dimensional. In particular, the degree map $\CH^2(X)\to\ZZ$ is far from being injective and the Beauville--Voisin subring $R(X)\subset \CH^*(X)$ is of infinite corank. 

However, for $k=\bar\QQ$ the situation should be drastically different. According to the general Bloch--Beilinson conjectures (see e.g.\ \cite{Bloch}) one expects that in this case $\CH^*(X)=R(X)$ (see \cite{HuyJEMS}).

Inspired by Bogomolov's `logical possibility' that every closed point of a K3 surface defined over $\bar\QQ$ may be
contained in a rational curve (see \cite{Bog}),
I have put forward (implicitly  in \cite{HuyJEMS, HuyMSRI} and explicitly in
related talks) the `categorial logical possibility' that for $k=\bar\QQ$ the derived category $\Db(X)$ might be generated by spherical objects. With our notation here, this means $\Db(X)=\ks^*$.

Bogomolov's logical possibility as well as $\Db(X)=\ks^*$ would clearly imply the Bloch--Beilinson conjecture
$\CH^2(X)\cong\ZZ$. To be precise, $\Db(X)=\ks^*$ would imply  $R(X)=\CH^*(X)$ which only
under the assumption $\rho(X)\geq2$ is known to yield $\CH^2(X)\cong\ZZ$.
To deduce  the Bloch--Beilinson conjecture from  Bogomolov's logical possibility one needs to use that for a point $x\in X$ contained
in a rational curve one has $[x]=c_X$ which was proved in \cite{BV}.

\medskip

The main purpose of this note is to show that the logical possibility $\Db(X)=\ks^*$ not only implies
the Bloch--Beilinson conjecture for K3 surfaces but that it is in fact equivalent to it. 
Note that there is no obvious relation between Bogomolov's geometric logical possibility and the categorial one.
Also, it is unclear whether Bogomolov's is actually equivalent to the Bloch--Beilinson conjecture.

\begin{thm}
Let $X$ be a projective K3 surface defined over $\bar \QQ$ with $\rho(X)\geq2$. Then the Bloch--Beilinson conjecture
for $X$, i.e.\ $\CH^2(X)\cong\ZZ$, holds if and only if $\Db(X)$ is generated by spherical objects, i.e.\
$\Db(X)=\ks^*$.
\end{thm} 

\begin{proof} 
If $\Db(X)=\ks^*$, then $\CH^*(X)=
K(X)=K(\Db(X))=K(\ks^*)$. By \cite{HuyJEMS} one has $K(\ks^*)=R(X)$ for $\rho(X)\geq2$ and
by \cite{BV} the cycle map is an isomorphism $R(X)\congpf\ZZ\oplus \NS(X)\oplus\ZZ$. Thus, $\CH^2(X)\cong \ZZ$.

The `only if' is more surprising. So suppose conversely that $\CH^2(X)\cong\ZZ$.
Then $R(X)=\CH^*(X)$. In other words, for any $E\in\Db(X)$ one has
$[E]\in R(X)=K(\ks^*)$ and hence $E\in\kd_{R(X)}=\ks^*$ by Corollary \ref{cor:sphcatThom}. Thus $\Db(X)=\ks^*$.
\end{proof}

\begin{remark}
i) For $\rho(X)=1$ the argument still shows that $\CH^2(X)\cong\ZZ$ implies $\Db(X)=\ks^*$,
but conversely $\Db(X)=\ks^*$ only yields $\CH^*(X)=K(\ks^*)$ which we do not control completely
for $\rho(X)=1$. 

ii) Note that for $X$ defined over an algebraically closed extension $\bar\QQ\subset k$ of positive transcendence degree
 $\CH^2(X)\ne\ZZ$ and hence $\Db(X)$ is not generated by spherical objects, i.e.\ $\Db(X)\ne\ks^*$.
\end{remark}

I am not aware of any technique that possibly could prove 
 $\Db(X)=\ks^*$ for $X$ over $\bar\QQ$.
However, as a consequence of Corollary \ref{cor:no-t}
 we can at least say the following.

\begin{cor}
Let $X$ be a projective K3 surface over $\bar\QQ$ with $\rho(X)\geq2$. Then $\CH^2(X)\cong\ZZ$
(Bloch--Beilinson) if and only if there exists a bounded t-structure on $\ks^*$.\qqed 
\end{cor}

\end{document}